\documentclass{ambp}
\usepackage{amssymb,latexsym, amscd}
\usepackage[all]{xy}
\usepackage{graphicx}

 \newtheorem{rem}[theorem]{Remark}

\usepackage{latexsym}
\usepackage{amssymb}
\newcommand{\ben}{\begin{equation}}
\newcommand{\een}{\end{equation}}

\newcommand{\integer}{\ensuremath{{\mathbb Z}}}
\newcommand{\naturals}{\ensuremath{{\mathbb N}}}
\newcommand{\real}{\ensuremath{{\mathbb R}}}
\newcommand{\complex}{\ensuremath{{\mathbb C}}}

\newcommand{\U}[1]{\ensuremath{{\mathrm U( #1 )}}}

\newcommand{\Aa}{{\mathcal A}}

\newcommand{\UU}{{\mathcal U}}

\newcommand{\SSS}{{\mathcal S}}

\newcommand{\LL}{\mathcal{L}}
\newcommand{\MM}{\mathcal{M}}
\newcommand{\HH}{\mathcal{H}}

\newcommand{\Proj}{\mathbf{P}}

\newcommand{\Xx}{\mathsf{X}}
\newcommand{\Mm}{\mathsf{M}}
\newcommand{\Gg}{\mathsf{G}}
\newcommand{\Hh}{\mathsf{H}}
\newcommand{\Ss}{\mathsf{S}}

\newcommand{\Kk}{\mathsf{K}}
\newcommand{\target}{\mathsf{t}}
\newcommand{\source}{\mathsf{s}}
\newcommand{\ident}{\mathsf{e}}

\newcommand{\mult}{\mathsf{m}}

\newcommand{\Loop}{\mathsf{L}}

\newcommand{\To}{\longrightarrow}

\newcommand{\timests}{\: {}_{\target}  \! \times_{\source}}

\newcommand{\toparrow}[1]{\stackrel{#1}{\longrightarrow}}

\newcommand{\VS}{\mathsf{Vector\ Spaces}}
\newcommand{\hol}{\mathrm{hol}}

\begin{document}



\title{An Introduction to  Gerbes on Orbifolds.}

\FirstAuthor{Ernesto Lupercio}
\FirstDepartment{Departamento de Matem\'{a}ticas}
\FirstInstitution{CINVESTAV}
\FirstStreetAddress{Apartado Postal 14-740 07000}
\FirstCity{M\'{e}xico D. F.}
\FirstCountry{M\'{e}xico}
\FirstSupport{The first author was partially supported by the National Science
Foundation and Conacyt-M\'exico}
\FirstEmail{lupercio@math.cinvestav.mx}
%
\SecondAuthor{Bernardo Uribe}
\SecondDepartment{Department of Mathematics}
\SecondInstitution{University of MIchigan}
\SecondStreetAddress{East Hall} 
\SecondCity{Ann Arbor, MI 48109}
\SecondCountry{USA}
\SecondEmail{uribeb@umich.edu}
\maketitle
\ShortTitle{Gerbes on Orbifolds}

\begin{abstract}
This paper is a gentle introduction to some recent results involving the theory of gerbes over orbifolds for topologists,
 geometers and physicists. We introduce gerbes on manifolds, orbifolds, the Dixmier-Douady class, Beilinson-Deligne orbifold
cohomology, Cheeger-Simons orbifold cohomology and string connections.
\end{abstract}

\tableofcontents

\section{Gerbes on smooth manifolds.}

We will start by explaining a well known example arising in electromagnetism as a motivation for the theory of gerbes. We will consider our space-time as canonically split as follows $$ M^4 = \real^ 4 = \real^3 \times \real = \{ (x_1, x_2, x_3 ; t) \colon x\in \real^3, t \in \real \}.$$ We will consider a collection of differential forms as follows
\begin{itemize}
   \item The electric field $E\in\Omega^1(\real^3)$.
   \item The magnetic field $B \in \Omega^2(\real^3)$.
   \item The electric current $J_E \in \Omega^2(\real^3)$.
   \item The electric charge density $\rho_E \in \Omega^3(\real^3)$.
\end{itemize}

We will assume that these differential forms {\emph{depend on $t$}} (so to be fair $ E \colon \real \to \Omega^1(\real^3)$, etc.).

We will {\emph{define}} the \emph{intensity of the electromagnetic field} by
$$ F = B - dt \wedge E \in \Omega^2(M)$$
and the compactly supported {\emph{electric current}} by
$$ j_E = \rho_E - dt \wedge J_E \in \Omega^3_c(M). $$
We are ready to write the {\textbf{Maxwell equations}}. They are
$$ dF = 0, \ \ \ \ \ \ \ \ \ \ \ \ \ d*F=j_E.$$
They are partial differential equations where the unknowns are the $3+3$ time-dependent components of the electric and the magnetic field.

If we would like them to look more symmetric we would need to introduce ``magnetic monopoles", namely a compactly supported 3-form for the magnetic charge density
$$j_B \in \Omega^3_c(M)$$
and rewrite the equations as
$$dF = j_B , \ \ \ \ \ \ \ \ \ \ \ \ d*F = j_E.$$

Now we let $N_t = \real^3 \times\{t\}$ be a space-like slice. Then the instantaneous total electric magnetic charges are respectively
$$ \int _{N_t}  j_E \ \ \ \ \ \ \  {\mathrm{and}} \ \ \ \ \ \  \int_{N_t} j_B.$$
But we prefer to consider the charges as elements in cohomology, namely
$$Q^t_E=[j_E|_{N_t}]\in H^3_c(N_t)$$
and
$$Q^t_B=[j_B|_{N_t}]\in H^3_c(N_t).$$

Now, quantum mechanics predicts that the charges above are quantized by the so-called \emph{Dirac quantization condition}, namely $Q_E^t$ is in the image of the homomorphism $$H_c^3(N_t,\integer) \to H_c^3(N_t;\real).$$ We can give a geometric interpretation to this quantization condition. For this purpose we must introduce the concept of (abelian) \emph{gauge field}.

\begin{definition} Let $M$ be a manifold. A $\U{1}$-gauge field on $M$ consists of a line bundle with a connection on $M$, namely
\begin{itemize}
\item[i)] A good Leray atlas $\UU= \{U_i\}_i$ of $M$.
\item[ii)] Smooth transition maps $g_{ij}\colon U_{ij}:=U_i \cap U_j \To \U{1}$. (These are the gluing maps that define the line bundle).
\item[iii)] A collection $(A_i)_i$ of 1-forms $A_i \in \Omega^1(U_i)$ that together are referred to as the field potential.
\item[iv)] These forms must satisfy the following equations:
  \begin{itemize}
      \item[a)] $g_{ij}$ is a cocycle (i.e. $g_{ij}g_{jk}=g_{ik}$ on $U_{ijk}$)
      \item[b)] $dA_i = dA_j$ on $U_{ij} = U_i \cap U_j$.
      \item[c)] $A_j-A_i = -\sqrt{-1} g_{ij}^{-1} d g_{ij}$.
   \end{itemize}
\item[v)] The 2-form $\omega=F=dA \in \Omega^2(M)$ is called the {\textbf{curvature}} of the connection $A$.
\end{itemize}
\end{definition}

It is an immediate consequence of the definition that the \emph{Bianchi identity} is satisfied, that is:
$$dF=0$$
and therefore we have a de Rham cohomology class $-[F] \in H^2(M,\real)$.

We can use the fact that $g_{ij}$ is a cocycle and consider its \v{C}ech cohomology class  $[g] \in H^1(M,{\underline{\U{1}}}) $
where ${\underline{\U{1}}}$ is considered as a sheaf over $M$. The exponential sequence of sheaves
$$ 0 \To \integer \To \real \toparrow{\exp( 2 \pi i \_)} \U{1} \To 1 $$
 immediately implies an isomorphism $$ H^1(M,{\underline{\U{1}}}) \cong H^2(M,\integer)$$
The class of $[g]$ in $H^2(M,\integer)$ is called \emph{the Chern class $c_1(L)$ of $L$}.

It is a theorem of Weil \cite{Weil52} that $-[F]$ is the image of the Chern class $c_1(L)$ under the map
$ H^2(M,\integer) \to H^2(M, \real)$. The Chern class completely determines the isomorphism type of the line bundle $L$, but does not determine the isomorphism class of the connection.

We say that a line bundle with connection is flat if its curvature vanishes. We have therefore that if a line bundle with connection is flat then its Chern class is a torsion class.

To solve the Maxwell equations is therefore equivalent to finding a line bundle with connection that in addition satisfies the field equation $d*F=j_E$. Let us for a moment consider the equation in the vacuum, namely consider the case of the field equation of the form $d*F=0$. We can write a rather elegant variational problem that solves the Maxwell equations in the vacuum (we learned this formulation from Dan Freed). Moreover, we can do so in a manner that exhibits fully the magnetic-electric duality of the problem. Let $A'$ be a second connection so that $F'=*F$. The electromagnetic Lagrangian is
$$ L(A,A') = \int_M \left(\frac{1}{4} |F|^2 + \frac{1}{4} |F'|^2 \right) dV$$
Clearly the equations in the vacuum are the Euler-Lagrange equation for $L(A,A')$, namely $\delta L=0$.

To add charges to the previous Lagrangian we consider a electrically charged particle whose worldline is a mapping $\gamma$ from a compact one-dimensional manifold to $M$. We consider the charge as an element $q\in H^0(\gamma, \integer) = \{ q | q\colon \gamma \To \integer \}$. To identify this with the charge as an element in $H^3_c(M,\integer)$ we us the Gysin map in cohomology $$i_! \colon H^0(\gamma,\integer) \To H^3_c(M,\integer)$$
given by the Thom-Pontrjagin collapse map and the Thom isomorphism. We can write the new Lagrangian that includes charges
$$ L = \int_M \left(\frac{1}{4} |B|^2 + \frac{1}{4} |B'|^2 \right) dV + i \int_\gamma \frac{1}{2} qA$$
Several remarks are in order.
\begin{itemize}
\item We have switched notations. We call $B$ what we used to call $F$. This is unfortunate but matches better the rest of the discussion.
\item It is no longer true that $dB=0$ (that is after all the whole point). In fact $B$ is no longer a global form.
\item Likewise $A$ is not a global form an actually only $\exp\left(i \int_\gamma qA\right)$ is well defined. Nevertheless the Lagrangian does define the correct Euler-Lagrange equations.
\end{itemize}

This situation is no longer a form of a line bundle with a connection. In spite of this, there is a geometric interpretation of the
 previous situation. This can be seen as a motivation for the introduction of the concept of \emph{gerbe} (cf. \cite{Hitchin}). (For more details on the physics see \cite{FreedHopkins, FreedDiff}.)

\begin{definition}\label{gerbemanifold} Let $M$ be a manifold. A {\textbf{gerbe}} with connection on $M$ is given by the following data:
\begin{itemize}
\item[i)] A good Leray atlas $\UU= \{U_i\}_i$ of $M$.
\item[ii)] Smooth maps $g_{ijk}\colon U_{ijk} \To \U{1}$.
\item[iii)] A collection $(A_{ij})$ of 1-forms $A_{ij} \in \Omega^1(U_{ij})$.
\item[iv)] A collection $B_i$ of 2-forms $B_i \in \Omega^2(U_i)$
\item[v)] These forms must satisfy the following equations:
  \begin{itemize}
     \item [a)]  $g_{ijk}$ is a cocycle (i.e. $g_{ijk}g_{ijl}^{-1}g_{ikl}g_{jkl}^{-1} =1$).
     \item [b)] $A_{ij}+A_{jk}-A_{ik} =-\sqrt{-1} d \log g_{ijk}$
     \item [c)] $B_j-B_i = dA_{ij}$
   \end{itemize}
\item[vi)] The global 3-form $\omega=dB \in \Omega^3(M)$ is called the {\textbf{curvature}} of the gerbe with connection $(g,A,B)$.
\end{itemize}
\end{definition}

The class $[g_{ijk} ] \in H^2(M,{\underline{\U{1}}}) \cong H^3(M,\integer)$  (where the isomorphism is induced by the exponential sequence of sheaves) is called the \emph{Dixmier-Douady} class of the gerbe and is denoted by $dd(g)$.  Just as before the class $[\omega]\in H^3(M,\real)$ in de Rham cohomology is the real image of the Dixmier-Douady class $dd(g)\in H^3(M,\integer)$.

Gerbes on $M$ are classified up to isomorphism by their Dixmier-Douady class $dd(g)\in H^3(M,\integer)$. This again ignores the connection altogether. In any case we have the following fact.

\begin{proposition}
An isomorphism class of a gerbe on $M$ is the same as an isomorphism class of an infinite-dimensional Hilbert projective bundle on $M$.
\end{proposition}
\begin{proof}
We will use Kuiper's theorem that states that the group $U(\HH)$ of unitary
operators in a Hilbert space $\HH$ is contractible, and therefore one has
$$\Proj(\complex^\infty)\simeq K(\integer,2)\simeq BU(1)\simeq U(\HH)/U(1)=\Proj U(\HH).$$ This fact immediately implies
$K(\integer,3)\simeq B\Proj U(\HH)$. Hence the class $dd(g) \in
H^3(X,\integer)=[X,K(\integer,3)]=[X,B\Proj U(\HH)]$ produces a Hilbert projective
bundle ${\bf E}$.
\end{proof}
In fact more is true. The collection of all gerbes in $M$ form a group under tensor product since $\U{1}$ is abelian
 (multiplication of the cocycles), and so do the set of all Hilbert projective bundles. One can prove that these two groups are isomorphic.

A gerbe with connection is said to be \emph{flat} if its curvature vanishes. Notice the following consequence of this fact,
\begin{proposition}
 A gerbe with connection is flat if and only if $dd(g)$ is a torsion class in cohomology.
 This is the case if and only if the projective bundle ${\bf E}$ is finite dimensional.
\end{proposition}
\begin{proof}
This is true because of a result of Serre \cite{DonovanKaroubi} valid for any CW-complex $M$. It states that
   if a class $\alpha \in H^3(M,\integer)$ is a torsion element
   then there exists a principal bundle $Z\to M$ with structure
   group $\Proj\U{n}$ so that when seen as an element
   $ \beta \in [M,B\Proj\U{n}] \to [M,B\Proj{\mathrm{U}}] = [M , B B \U{1}] = [M , B
   K(\integer,2)] = [M, K(\integer, 3)]= H^3(M,\integer)$ then
   $\alpha = \beta$. In other words, the image of $[M,B\Proj\U{n}]
   \to H^3(M,\integer)$ is exactly the subgroup of torsion elements that are killed by multiplication by $n$.
\end{proof}

We refer the reader to the papers \cite{Murraygerbes, Murray02} for gerbes from the point of view of bundle gerbes.

\section{Orbifolds}

The notion of orbifold was first introduced by Satake in his
seminal paper \cite{Satake56}. In this 1956 paper Satake defines
for the very first time the concept of an orbifold by means of
orbifold atlases whose charts Satake calls local uniformizing
systems. The name that orbifolds take in this early work are
\emph{$V-$manifolds}. Quite remarkably he already works with a
version of \v Cech groups. He goes on to prove the De Rham theorem
and Poincar\'e duality with rational coefficients. For about two
decades the japanese school carried out brilliantly the study of
orbifolds. It deserves special mention the work of Tetsuro
Kawasaki. In his papers of the late 70's Kawasaki generalizes
index theory to the orbifold setting \cite{Kawasaki78, Kawasaki79,
Kawasaki81}. Another important work along these veins in the work
of Thurston specially his concept of orbifold fundamental group
\cite{Thurston97}

Somewhat independently the algebraic geometers developed the
concept of \emph{stack} in order to deal with moduli problems. As
it happens orbifolds arise quite naturally from the very same
moduli problems and it didn't take long to realize that the theory
of stacks provided another way of understanding the category of
orbifolds, and viceversa. For example, the Deligne-Mumford moduli
stack $\MM_g$ for genus $g$ curves \cite{DeligneMumford} is in
fact an orbifold. This is one of the reasons for the importance of
orbifolds, many moduli spaces are better understood as orbifolds.
The paper of Artin \cite{Artin74} is the place where a very
explicit conection with groupoid atlases takes place
for the first time. Implicitly these ideas are already  present in
Grothendieck's toposes \cite{Grothendieck72}. The groupoiod
approach to orbifolds is finally carried out by Haefliger
\cite{Haefliger} and by Moerdijk and his collaborators
\cite{Moerdijk91,MoerdijkPronk1,CrainicMoerdijk,MoerdijkPronk}. In
this work they put forward the important concept of Morita
equivalence.

 The interest of orbifolds in physics can be traced back to the work of Dixon, Harvey, Vafa and
Witten \cite{DHVW,DHVW2} where motivated by superstring
compactification they introduce a orbifold theory using a K3 with
27 singular points. It is there that the orbifold Euler
characteristic is defined motivated by the physics. It is a
remarkably insightful notion of their work to realize that their
results depend only on the orbifold and not on group actions, for
all their example are global orbifolds.  This work produced an
explosion of activity related to orbifolds in the physics
community. The introduction to the mathematics side of the
geometrization of many of these ideas and results is due to Chen
and Ruan. Their highly influential papers \cite{ChenRuan,Ruan}
introduced many concepts from the physics literature rigorously
into symplectic and algebraic geometry. In this work orbifolds are
completely general, not necessarily global quotients. In
particular they discovered a remarkable cohomology (the Chen-Ruan
cohomology) that was never looked for by mathematicians before,
that is amenable to mirror symmetry and is the object of very
intense research as we write this. Ruan himself introduced
twistings (related to the $B$-field in string theory) into his
theory, and it was his work was motivated us to consider gerbes
over orbifolds in the first place.

The definition of an orbifold is quite involved. We will start by trying to motivate the definition with a few examples.

\begin{example}\label{pillow}
Let $M=T^2 = S^1 \times S^1$ be a two-dimensional torus, and let $G= \integer_2$ be the \emph{finite} subgroup of diffeomorphisms of $M$ given by the action
$$ (z,w) \mapsto (\bar{z}, \bar{w})$$
It is not hard to show that while the quotient space $X=M/G$ is topologically a sphere it is impossible to put a smooth structure on $X$ so that the quotient map $M \To X$ will become smooth. It is in this sense that we say that $X$ is not a smooth manifold.

What we can still do is to enlarge the category of smooth manifolds to a bigger category called the category of \emph{orbifolds}. Once we do this by considering the orbifolds $\Mm$ and $\Xx$ the natural orbifold map $\Mm \To \Xx$ is smooth.

\begin{eqnarray*}
\includegraphics[height=1.5in]{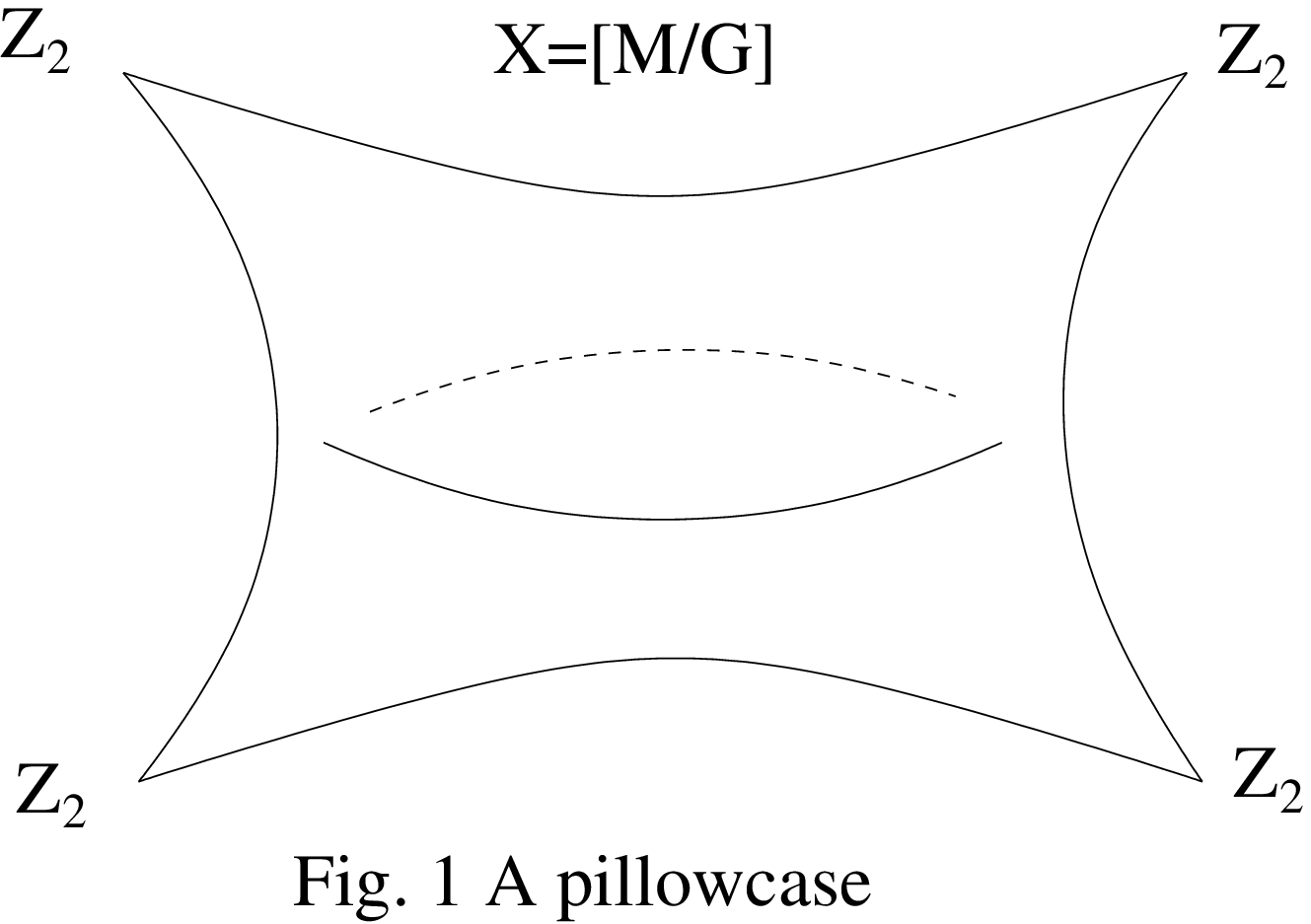}
\end{eqnarray*}

While the orbifold $\Mm$ contains exactly the same amount of information as $M$ the orbifold $\Xx=[M/G]$ (known as a pillowcase) contains more information that the quotient space $X=M/G$. For example $\Xx$ remembers that the action had 4 fixed points each with stabilizer $G$. It remembers in fact the stabilizer of every point, and how these stabilizers fit together. On the other hand $\Xx$ does not remember neither the manifold $M$ nor the group $G$. In fact if we define $N$ to be two disjoint copies of $M$ and $H = G \times G$ to act on $M$ by letting $G\times1$ act by conjugation on both copies as before, and $1\times G$ act by swapping the copies then $$\Xx = [M/G] = [N/H].$$
\end{example}

\begin{example}
Not every orbifold can be obtained from a finite group acting on a manifold. An orbifold is always \emph{locally} the quotient of a manifold by a finite group but this may fail globally. For example consider the teardrop:

\begin{eqnarray*}
\includegraphics[height=1.5in]{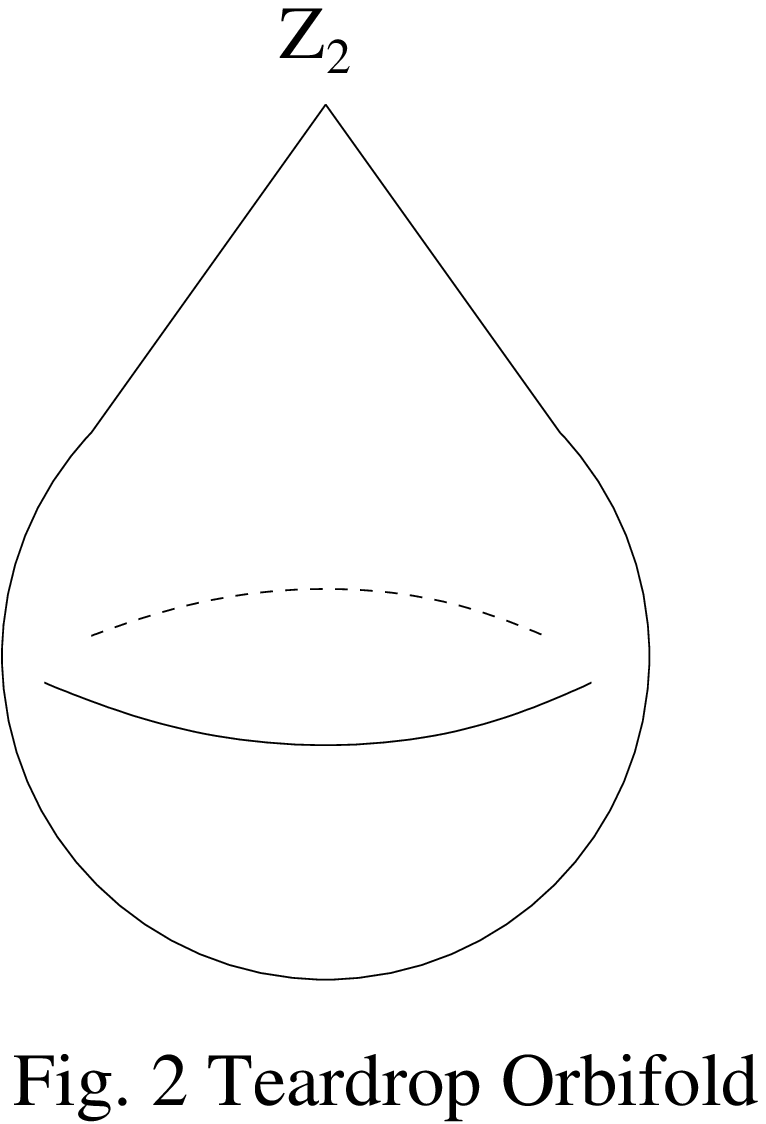}
\end{eqnarray*}

This orbifold may be obtained by gluing two global quotients. Consider the orbifold
 $\Xx_1 = [\complex / \integer_2]$ where $\integer_2$ acts by the holomorphic automorphism
 $z \mapsto -z$. Let $\Xx_2 = \complex$ simply be the complex plane.
Then we have in the category of orbifolds a diagram of inclusions
$$ \Xx_1 \longleftarrow \complex^* \longrightarrow \Xx_2 $$
and therefore we can glue $\Xx_1$ and $\Xx_2$ along $\complex^*$ obtaining the teardrop $\Xx$.

Notice that by the Lefschetz fixed point theorem it is not possible to obtain the teardrop as a global quotient of a manifold by a finite group.
Nevertheless it is possible to recover the teardrop as the quotient of the manifold $\complex^2 \backslash \{0\}$ by the action of the lie group
$\complex^*$. The action being $(z_1,z_2) \cdot \lambda \mapsto (\lambda z_1, \lambda^2z_2)$.
\end{example}

\begin{example} The present example only deals with smooth manifolds, so let $M$ be a smooth manifold.
 It is well known that a smooth manifold is a pair $(M,\UU)$ of a (Hausdorff, paracompact) topological
 space $M$ together with an atlas $\UU=\{U_i\}_{i\in I}$, and is only by abuse of notation that we speak of
 a manifold $M$. In fact a smooth manifold is actually an equivalence class of a pair $ [M,\UU]$ where we say that
 $(M,\UU_1) \sim (M,\UU_2)$ if and only if there is a common refinement $(M,\UU_3)$ of the atlas.

We can say this in a slightly different way that will be easier to generalize to the case of orbifolds. To have a pair $(M,\UU)$ is the same thing as to have a small topological category $\Mm_\UU$ defined as follows.
\begin{itemize}
   \item Objects: Pairs $(m,i)$ so that $m \in U_i$. We endow the space of objects with the topology $$ \coprod_i U_i.$$
   \item Arrows: Triples $(m,i,j)$ so that $m \in U_i \cap U_j =U_{ij}$. An arrow acts according to the following diagram. $$(x,i)\toparrow{(x,i,j)}(x,j).$$ The composition of arrows is given by $$(x,i,j)\circ(x,j,k)=(x,i,k)$$ The topology of the space of arrows in this case is $$ \coprod_{(i,j)} U_{ij}.$$

\end{itemize}

The category $\Mm$ is actually a {\textbf{groupoid}}, namely every arrow has an inverse, in fact $$(x,i,j) \circ (x,j,i) =(x,i,i) = Id_{(x,i)}.$$

We will therefore define a manifold to be the equivalence class of the groupoid $\Mm_\UU$ by a relation called Morita equivalence (that will amount exactly to the equivalence of atlases). We will define this equivalence relation presently.

\end{example}

As we said before a groupoid $\Gg$ in great generality is a category in which every
morphism is invertible. We will assume that $\Gg_0$ and $\Gg_1$, the sets of
objects and morphism respectively, are \emph{smooth manifolds}. We will denote the structure maps by:
      $$\xymatrix{
        \Gg_1 \timests \Gg_1 \ar[r]^{m} & \Gg_1 \ar[r]^i &
        \Gg_1 \ar@<.5ex>[r]^s \ar@<-.5ex>[r]_t & \Gg_0 \ar[r]^e & \Gg_1
      }$$
where $s$ and $t$ are the source and the target maps, $m$ is the composition (we can compose two arrows whenever the
target of the first equals the source of the second), $i$ gives us
the inverse arrow, and $e$ assigns the identity arrow to every
object. We will assume that all the structure maps are smooth maps. We also require that the
maps $s$ and $t$ must be submersions, so that $\Gg_1 \timests
\Gg_1 $ is also a manifold.

A topological (smooth) groupoid is called {\it \'{e}tale} if the
source and target maps $s$ and $t$ are local homeomorphisms (local
diffeomorphisms). For an \'{e}tale groupoid we will mean a
topological \'{e}tale groupoid.  We will always denote groupoids
by letters of the type $\Gg,\Hh,\Ss$. We will also assume that the  anchor map $(s,t):
\Gg_1 \to \Gg_0\times \Gg_0$ is proper, groupoids with this
property are called {\it proper groupoids}.  A theorem of Moerdijk and
Pronk \cite{MoerdijkPronk}
 states that the category of orbifolds is equivalent to a quotient category of the category of proper \'{e}tale groupoids
after inverting Morita equivalence. Whenever we write orbifold, we will choose a proper
\'{e}tale smooth groupoid representing it (up to Morita
equivalence.)

\begin{definition} A morphism of groupoids $\Psi: \Hh \to \Gg$ is a pair of maps
$\Psi_i: \Hh_i \to \Gg_i$ $i=0,1$ such that they commute with the
structure maps. The maps $\Psi_i$ will be required to be smooth.

The morphism $\Psi$ is called {\it Morita} if the following square
is a cartesian square
\begin{eqnarray}
 \xymatrix{
         \Hh_1 \ar[r]^{\Psi_1} \ar[d]_{(s,t)} & \Gg_1 \ar[d]^{(s,t)} \\
         \Hh_0 \times \Hh_0 \ar[r]^{\Psi_0 \times \Psi_0} & \Gg_0 \times \Gg_0
         } \label{Moritasquare}
\end{eqnarray}
and the map $s\pi_2 \colon \Hh_0{}_{\Psi_0}\times_t \Gg_1 \to \Gg_0$ to be
surjective and \'etale (local diffeomorphism).

 Two groupoids $\Gg$ and $\Hh$ are Morita equivalent if there
exist another groupoid $\Kk$ with Morita morphisms $\Gg
\stackrel{\simeq}{\leftarrow} \Kk \stackrel{\simeq}{\to} \Hh$.
\end{definition}

We often need a particular kind of representative for an orbifold,

\begin{definition}A groupoid $\Gg$ is called {\emph{ Leray}} if $\Gg_i$ is
diffeomorphic to a disjoint union of contractible open sets for
all $i \in \naturals$ where $\Gg_i = \Gg_1 \timests \Gg_1 \timests \cdots \timests \Gg_1 = \{ (\alpha_1,\ldots,\alpha_i) \colon s(\alpha_{j+1})=t(\alpha_i)\}$.\end{definition}
The existence of such Leray groupoid representative for every
orbifold is proved by Moerdijk and Pronk  \cite[Cor.
1.2.5]{MoerdijkPronk1}.

\begin{example}
   Consider again example \ref{pillow}. Define the following groupoids.

   \begin{itemize}
      \item The groupoid $\Gg$ whose space of objects are elements $m \in M$ with the topology of $M$, and whose space of arrows is the set of pairs $(m,g)$ with the topology of $M\times G$. We have the diagrams $$m\toparrow{(m,g)} mg$$ and the composition law $$(m,g) \circ (mg, h) = (m,gh).$$

      \item Similarly we define the groupoid $\Hh$ using the action of $H$ in $N$ with objects $n\in N$ and arrows $(n,h)\in N\times H$.
   \end{itemize}

The orbifold $\Xx$ is the equivalence class of the groupoid $\Gg$. Since $\Gg$ and $\Hh$ are Morita equivalent we can say equivalently that $\Xx$ is the equivalence class of $\Hh$. By abuse of notation we will often say that $\Gg$ is an orbifold when we really mean that its equivalence class is the orbifold.

\end{example}

\begin{example}\label{LerayGlobal}

   More generally, let $M$ be a smooth manifold and $G \subset {\mathrm{Diff}}(M)$ be a finite group acting on it.

\begin{itemize}

\item  We say that the orbifold $[M/G]$ is the equivalence class of the groupoid $\Xx$ with objects $m \in M$ and arrows $(m,g) \in M\times G$.

\item  We can define another groupoid representing the same orbifold as follows. Take a contractible open cover $\UU=\{U_i\}_{i \in I}$
of $M$ such that all the finite intersections of the cover are either contractible or empty, and with
the property that for any $g \in G$ and any $i \in I$ there exists $j \in I$ so that $U_ig=U_j$. Define
$\Gg_0$ as the disjoint union of the $U_i$'s with $\Gg_0 \stackrel{\rho}{\to} M=\Xx_0$ the natural map.
Take $\Gg_1$ as the pullback square
  $$\xymatrix{ \Gg_1 \ar[r] \ar[d] & M \times G \ar[d]^{\source \times \target} \\
             \Gg_0 \times \Gg_0 \ar[r]^{\rho \times \rho}& M \times M}$$
where $\source(m,g) =m$ and $\target(m,g) = mg$.
From the construction of $\Gg$ we see that we can think of $\Gg_1$ as the disjoint union of all the intersections
of two sets on the base times the group $G$, i.e.
$$\Gg_1 = \left( \bigsqcup_{(i,j) \in I\times I} U_i \cap U_j \right) \times G$$
where the arrows in $U_i \cap U_j \times\{g\}$ start in $U_i|_{U_j}$ and end in $(U_j|_{U_i})g$.  This defines the proper \'{e}tale Leray groupoid $\Gg$ and
by definition it is Morita equivalent to $\Xx$.
\end{itemize}

\end{example}

Given an orbifold a very important construction is that of its \emph{classifying space}.  The nerve of a groupoid (see \cite{Segal1})  is a semisimplicial set $N\Gg$ where the objects of $\Gg$ are
the vertices, the morphisms the 1-simplexes, the triangular commutative diagram the 2-simplexes, and so on.
We can define the boundary maps $\delta_i : \Gg_i \to \Gg_{i-1}$ by:

\begin{displaymath}
\delta_i(x_1, \dots , x_n) = \left\{
 \begin{array}{ll}
(x_2, \dots , x_n) & \mbox{if $ i=0$} \\
(x_1, \dots, m(x_i,x_{i+1}), \dots , x_n) & \mbox{ if $ 1 \leq i \leq n-1$} \\
(x_1, \dots, x_{n-1}) & \mbox{if $ i =n $} \\
\end{array}
 \right\}
\end{displaymath}

$N \Gg$ determines $\Gg$ and its geometric realization $B \Gg$ is called the classifying space of the orbifold. This space is important to us because it is a result of Moerdijk \cite{Moerdijk98} that $H^*(\Xx; \integer)
\cong H^*(B\Xx; \integer)$ where the left hand side is sheaf cohomology (to be defined) and the right hand side is simplicial cohomology of
$B\Xx \simeq M_G$ with coefficients in $\integer$.

\begin{example} In the case in which $\Gg=[M/G]$ then $B\Gg = M \times_G EG$ is the Borel construction associated to the group action.

A particular but nevertheless important example is that of $\Gg=[*/G]$ a finite group acting on a point.
 In this case $B \Gg = BG$. The result of Moerdijk mentioned above implies that $BG$ computes group cohomology.

\end{example}

\begin{example} This example is due to Segal \cite{Segal1}. If $(M,\UU)$ is an atlas for $M$ then the classifying space of the groupoid $\Mm_\UU$ is $$ B \Mm_\UU \simeq M.$$\end{example}

\section{Gerbes over Orbifolds.}

In this section we discuss definitions and result first introduced in \cite{LupercioUribeKtheory}.

\begin{example} Let us recast the definition of gerbe over a manifold $(M,\UU)$, with Leray groupoid $\Mm_\UU$. Notice than in this case
\begin{itemize}
   \item $(\Mm_\UU)_0 = \coprod_i U_i$
   \item $(\Mm_\UU)_1=\coprod_{(i,j)} U_{ij}$
   \item $(\Mm_\UU)_2=\coprod_{(i,j,k)} U_{ijk}$
\end{itemize}
and so on.

 To have a gerbe over an orbifold is the same as to have a map $g:(\Mm_\UU)_2 \To \U{1}$ satisfying the cocycle condition. The data defining a gerbe with connection are in addition forms $A\in \Omega^1((\Mm_\UU)_1)$ and $B\in \Omega^2((\Mm_\UU)_0$, satisfying the equations of definition \ref{gerbemanifold}
\end{example}

\begin{definition} A gerbe (with band $\U{1}$) over an orbifold is a pair $(\Gg,g)$ where $\Gg$ is a groupoid representing the orbifold and $g$ is
a 2-cocycle $g \colon \Gg_2 \to \U{1}$. A gerbe with connection consists of
a 1-form $A \in \Omega^1(\Gg_1)$, a 2-form $B \in
\Omega^2(\Gg_0)$  satisfying:
\begin{itemize}
\item $\target^*B - \source^*B = dA$ and
\item $\pi_1^*A + \pi_2^*A -\mult^*A = -\sqrt{-1} g^{-1} dg$
\end{itemize}
 The $\Gg$-invariant 3-form $\omega=dB \in \Omega^3(\Gg_0)$ is called the curvature  of the gerbe with connection $(g,A,B)$. Here by $\Gg$-invariant we mean that $\source^* \omega = \target^* \omega$.
\end{definition}

The following theorem of \cite{LupercioUribeKtheory, LupercioUribeHolonomy} describes the basic classification of gerbes over orbifold (without a connection).

\begin{theorem} The following holds.
\begin{itemize}
   \item Every gerbe on an orbifold has a representative of the form $(\Gg, g)$ where $\Gg$ is a Leray groupoid.
   \item We define the characteristic class $\ell(g)$ of $g$ to be the class in $H^3(B\Gg,\integer) \simeq H^3(\Gg,\integer)
 \simeq H^2(\Gg,{\underline{\U{1}}})$ induced by the \v{C}ech cocycle  $g \in C^2(\Gg,\underline{\U{1}})$.
 Then  isomorphism classes of gerbes  over the orbifold $\Gg$ are in one to one correspondence with $H^3(B\Gg,\integer)$ via the class $\ell(g)$.
\end{itemize}
\end{theorem}

To classify gerbes with connection $(g,A,B)$ up to isomorphism we need to introduce a new type of cohomology. We define now the so-called \emph{Beilinson-Deligne cohomology} of $\Gg$.

For the purpose of exposition we will introduce this cohomology for the Leray groupoid of Example \ref{LerayGlobal} and refer the reader to \cite{LupercioUribeBfield, LupercioUribeHolonomy} for the case of a general orbifold groupoid.

A $\Gg$-sheaf is a sheaf over $\Gg$ on which $\Gg$ acts continuously. Let $\Aa^p_{\Gg}$ denote the $\Gg$-sheaf of differential $p$-forms and $\integer_{\Gg}$ the constant $\integer$ valued $\Gg$
sheaf with $\integer_{\Gg} \to \Aa^0_{\Gg}$ the natural inclusion of constant into smooth functions.

Let's denote by
$\breve{C}^*(\Gg;\U{1}(q))$ the total complex
\begin{eqnarray*}
\xymatrix{
\breve{C}^0(\Gg; \U{1}(q))  \ar[r]^{\delta -d } &
         \breve{C}^1(\Gg;\U{1}(q) )
        \ar[r]^{\delta + d}  &  \breve{C}^2(\Gg; \U{1}(q))
         \ar[r]^{\ \ \delta-d} & \cdots
}
\end{eqnarray*}
induced
by the double complex
\begin{eqnarray} \label{double complex}
   \xymatrix{
     \vdots  & & \vdots & & \vdots\\
     \Gamma(\Gg_2, \U{1}_\Gg) \ar[u]^\delta \ar[rr]^{-\sqrt{-1}d \log} & &
         \Gamma(\Gg_2, \Aa^1_{\Gg})
        \ar[r]^d \ar[u]^\delta  & \cdots \ar[r]^d & \Gamma(\Gg_2, \Aa^{q-1}_{\Gg})
         \ar[u]^\delta\\
      \Gamma(\Gg_1, \U{1}_\Gg) \ar[u]^\delta \ar[rr]^{-\sqrt{-1}d \log} & &
          \Gamma(\Gg_1, \Aa^1_{\Gg})
        \ar[r]^d \ar[u]^\delta  & \cdots \ar[r]^d & \Gamma(\Gg_1, \Aa^{q-1}_{\Gg})
        \ar[u]^\delta\\
      \Gamma(\Gg_0, \U{1}_\Gg) \ar[u]^\delta \ar[rr]^{-\sqrt{-1}d \log} & &
         \Gamma(\Gg_0, \Aa^1_{\Gg})
          \ar[r]^d \ar[u]^\delta  & \cdots\ar[r]^d & \Gamma(\Gg_0, \Aa^{q-1}_{\Gg})
          \ar[u]^\delta
    }
\end{eqnarray}
with $(\delta +(-1)^{i} d )$ as coboundary operator, where the $\delta$'s are the maps induced by the
simplicial structure of the nerve of the category $\Gg$ and $\Gamma(\Gg_i, \Aa^j_{\Gg})$ stands
for the global sections of the sheaf that induces $\Aa^j_{\Gg}$ over $\Gg_i$ (see \cite{LupercioUribeHolonomy}).
Then the Beilinson-Deligne cohomology is defined as as follows:
$$H^n(\Gg, \integer(q)) \cong H^{n-1}(\Gg, \U{1}(q)) :=H^{n-1}\breve{C}(\Gg;\U{1}(q)).$$

It is proved in \cite{LupercioUribeBfield} that $H^n(\Gg,\integer(q))$ only depends on the orbifold and not on the particular groupoid used to represent it. In the same paper the notation $H^n(\Gg,\integer(q))$ (given by a refined version of the exponential sequence of sheaves for complexes of sheaves) is explained.

We have the following.

\begin{proposition}
For $\Gg$ a Leray description of a smooth \'{e}tale groupoid, a gerbe with
connection is a 2-cocycle of the complex $\breve{C}(\Gg,\U{1}(3))$, that
is, a triple $(h,A,B)$ with $B \in \Gamma(\Gg_0, \Aa^2_{\Gg})$,
$A \in \Gamma(\Gg_1, \Aa^1_{\Gg})$ and
$h \in \Gamma(\Gg_2, \U{1}_\Gg)$ that satisfies $\delta B = dA$,
$\delta A = -\sqrt{-1}d \log h$ and $\delta h = 1$.
\end{proposition}

\begin{definition}
An $n$-gerbe with connective structure over $\Gg$ is an
$(n+1)$-cocycle of $\breve{C}^{n+1}(\Gg,\U{1}(n+2))$.
Their isomorphism classes are classified by
$$H^{n+1}(\Gg,\U{1}(n+2)) =
H^{n+2}(\Gg,\integer(n+2)).$$
\end{definition}

The following theorems were proved in \cite{LupercioUribeBfield, LupercioUribeHolonomy}.

\begin{proposition} \label{uninteresting part}
$$H^{p}(\Gg,\integer(n)) \cong H^{p-1}(\Gg,\U{1}(n)) = \left\{
\begin{array}{cc}
H^{p-1}(\Gg, \underline{\U{1}}) = H^{p}(\Gg, \integer) & \mbox{for} \ p > n \\
H^{p-1}(\Gg, \U{1}) & \mbox{for} \ p < n
\end{array} \right.$$
where $\underline{\U{1}}$ stands for the sheaf of $\U{1}$ valued functions.
\end{proposition}

We have argued in \cite{LupercioUribeBfield} that a $B$-field in the physics terminology for type II orbifold superstring
theories is the same as a gerbe with connection on the orbifold.

The following theorem generalizes a result of Brylinski that he proved in the case of a smooth manifold $M$ \cite{Brylinski}.

\begin{theorem} We have the following classifications.
\begin{itemize}
    \item The group of isomorphism classes of line orbibundles with connection on $\Gg$ is isomorphic to $H^2(M,\integer(2))$.
   \item The group of isomorphism classes of gerbes with connection on $\Gg$ is isomorphic to $H^3(M,\integer(3))$.
\end{itemize}
\end{theorem}

\begin{rem} It is quite interesting to point out that if $[g,A,B]$ is the BD-class of $(g,A,B)$ then $\omega = dB$ is completely determined by $[g,A,B]$. We call the 3-from $\omega$ the curvature of the class $[g,A,B]$. An analogous definition can be made for $n$-gerbes yielding a $(n+2)$-form $\omega$. \end{rem}

A \emph{discrete torsion} on an orbifold $\Gg=[M/G]$ is a 2-cocycle $\theta \colon G\times G \to U(1)$ in the bar group cohomology complex of $G$ \cite{VafaWitten} (cf. \cite{Sharpe}).

\begin{proposition}\cite{LupercioUribeBfield} For a global orbifold $[M/G]$ the map $\theta \mapsto (\theta,0,0)$ injects the group of discrete torsions of an orbifold into the group of flat gerbes (=flat B-fields). In fact the induced map in cohomology $H^3(G,\integer) \To H^3(\Gg,\integer(3))$ is injective. \end{proposition}

\begin{rem}
Let us remark that the gerbes coming from discrete torsion do not amount to all the flat gerbes.
Consider the case in which $G=\{1\}$ and $H^2(M, \U{1}) \neq 0$, then there is no discrete torsion
but there are non trivial flat gerbes.
\end{rem}

\section{Holonomy}

To warm up consider a line bundle with connection $(L,g,A)$ over a manifold $(M,\UU)$. Classically the holonomy of $(L,g,A)$ determines for every path $\gamma \colon [0,T] \To M$ a linear mapping $$\hol_{(L,g,A)}(\gamma) \colon L_{\gamma(0)} \To L_{\gamma(T)}$$ that composes well with path concatenation. On a chart $\gamma \colon [0,T] \To V\in \real^n$ of $M$ where $L=V\times \complex$ we can write such a map simply  as an element in $\U{1}$ by
$$ \hol_{(L,g,A)}(\gamma) = \exp \left( 2 \pi i \int_\gamma A \right). $$
This formula is enough to completely define the holonomy for manifolds in general in view of the following.
\begin{proposition} Let $\SSS^0(M)$ be the 0-th Segal category of $M$ having
\begin{itemize}
\item Objects: The points $m\in M$.
\item Arrows: Paths $\gamma \colon [0,T] \To M$ with composition given by concatenation of paths.
\end{itemize}
Then the holonomy of a line bundle with connection defines a functor $$\hol_{(g,A)} \colon \SSS^0(M) \to \VS_1(\complex)$$
from $\SSS^0(M)$ to the category of 1-dimensional vector spaces with linear isomorphisms.
\end{proposition}

Notice that we can restrict our attention to the closed paths (automorphisms of $\SSS^0(M)$) to obtain a function on the loop space $\LL M$
of $M$
$$ \hol^\circ_{(g,A)} \colon \LL M \to \U{1}$$
We consider this function as an element $\hol^\circ_{(g,A)} \in H^0(\LL M, {\underline{\U{1}}})$.

\begin{definition} The \emph{transgression map}
 $H^2(M; \integer) \rightarrow H^1(\LL M ; \integer)$
is defined as the following composition. Let $$S^1 \times \LL M \longrightarrow M$$ be the
\emph{evaluation map} sending $(z,\gamma) \mapsto \gamma(z)$. We
can use this map together with the K\"unneth theorem and the fact
that $H^1(S^1;\integer)=\integer$  to get
$$H^2(M; \integer) \rightarrow H^2(S^1 \times \LL M ; \integer) \cong
H^2(\LL M; \integer) \oplus ( H^1(\LL M ; \integer) \otimes
H^1(S^1; \integer) )$$ $$ \stackrel{\cong}{\rightarrow} H^2(\LL M;
\integer) \oplus H^1(\LL M ; \integer) \rightarrow H^1(\LL M ;
\integer)\cong H^0(\LL M ; \underline{\U{1}})$$ (where
the next to last map is projection into the second component, and
the last is induced by the exponential sequence). \end{definition}

\begin{proposition} The element
$\hol^\circ_{(g,A)} \in H^0(\LL M, {\underline{\U{1}}})$
is the image of $c_1(g)\in H^2(M,\integer)$ under the transgression map.
\end{proposition}

This implies that  $\hol^\circ_{(g,A)}$ depends only on the Chern class (namely on the isomorphism class of $(L,g)$ and not on the specific connection $A$. So the functor $\hol_{(g,A)}$ contains more information that $\hol^\circ_{(g,A)}$.

\begin{example} Suppose that $\omega = dA=0$, so the line bundle $L$ is flat. Then $c_1(g)$ is a torsion class. In this case the holonomy induces a homomorphism $\rho \colon \pi_1(M) \To \U{1}$ that determines the functor $\hol_{(g,A)}$ up to natural transformation. \end{example}

 Let us consider consider the holonomy as a map
$$ \hol^Z_{(g,A)} \colon Z_1(M) \To \U{1},$$
where $Z_1(M)$ are the closed smooth 1-chains on $M$.
We define  $\chi$ to be $$\chi:= -\frac{\sqrt{-1}}{2 \pi}\log \hol^Z.$$ If we consider the curvature of $L$ as
 a 2-form $\omega$ on $M$ we have obtained a pair $(\chi, \omega)$ with
 $$\chi \colon Z_1 (M) \To \real / \integer$$
 and
 $$\chi(\partial c ) =  \int_c \omega \hbox{ \ mod \ } \integer$$
 whenever $c$ is a smooth 2-chain (the pair $(\chi, \omega)$ is called a differential character).

Following Cheeger-Simons \cite{CheegerSimons} we will denote by $\hat{H}^2 _{cs}(M)$
 the group of such differential characters of $M$.

If we substitute
 the line bundle by a $(q-2)$-gerbe with connection. The holonomy becomes now a homomorphism $Z_{q-1}(M) \to \U{1}$, then we can define in general  $\hat{H}^q _{cs}(M).$

The following theorem \cite{Brylinski, Murraygerbes} relates the CS-cohomology  to the BD-cohomology of a manifold $M$:
\begin{theorem}
$$H^{q}(M; \integer(q)) \cong  \hat{H}^q _{cs}(M).$$
\end{theorem}

Actually the holonomy of a gerbe can also be seen as a functor.

\begin{theorem}  Let $\SSS^1(M)$ be the 1-st Segal category of $M$ having
\begin{itemize}
\item Objects: Maps $\gamma \colon S^1 \coprod \ldots \coprod S^1 \To M$.
\item Arrows: Maps $\Sigma \colon F \To M$ from 2-dimensional compact manifolds $F$ to $M$ forming cobordisms
   between two objects, with composition given by concatenation of surfaces.
\end{itemize}
Then the holonomy of a gerbe with connection $(g,A,B)$ defines a functor $$\hol_{(g,A,B)} \colon \SSS^1(M) \to \VS_1(\complex)$$
from $\SSS^1(M)$ to the category of 1-dimensional vector spaces with linear isomorphisms. Such a functor is called a \emph{string connection}.
\end{theorem}

For instance, in
   the picture below we have four maps $\gamma_i \colon S^1
   \rightarrow M$ ($i=1,2$) and a map $\Sigma \colon F \rightarrow
   M$ from a 2-dimensional manifold $F$ into $M$. Such a
   configuration would produce a linear isomorphism
   $$\hol_{(g,A,B)}(\Sigma) \colon L_{\gamma_1}\otimes L_{\gamma_2}  \longrightarrow
   L_{\gamma_3}\otimes L_{\gamma_4}.$$ Where $L$ is a line bundle on $\LL M$ defined by the functor. The reader may imagine that these are two strings evolving and interacting in $M$ if she prefers to do so.

\begin{eqnarray}
\includegraphics[height=1.0in]{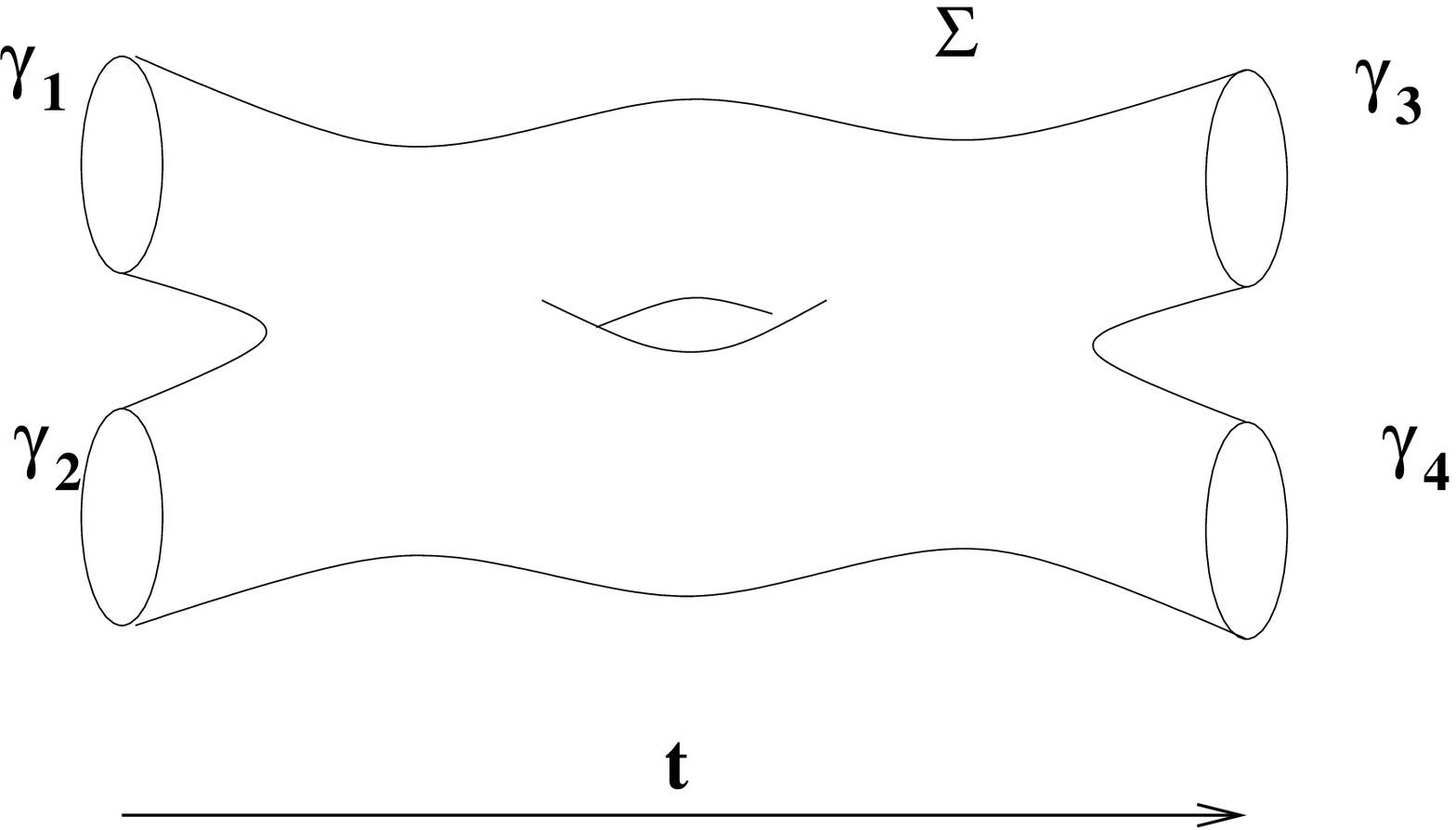} \label{graph string}
\end{eqnarray}

Now consider in general an orbifold $\Xx$. We will describe now the results of \cite{LupercioUribeLoop, LupercioUribeHolonomy, LupercioUribeStringI, LupercioUribeStringII} that refine the previous results to the case of orbifolds.

First we have defined an infinite dimensional orbifold, \emph{the loop orbifold} $\Loop \Xx$ associated to $\Xx$ by giving an explicit groupoid representation of it that we call the loop groupoid. This is quite technical in the general case \cite{LupercioUribeLoop}. In the global quotient case $[M/G]$ the situation is simpler. Let $\Gamma$ be a finite group (we may need to suppose $M$ is compact).

\begin{definition}
An orbifold loop $[M/G]$ will consist of a map $\phi \colon Q \to M$ of a $\Gamma$-principal
bundle $Q$ over the circle $S^1$ together with a homomorphism $\phi_\# \colon \Gamma \to G$ such that
$\phi$ is $\phi_\#$-equivariant. Let's denote this space of orbifold loops $(\phi,\phi_\#)$ by $\LL [M/G]$. It has a natural
action of the group $G$ as follows. For $h \in G$ let $\psi := \phi \cdot h$ where
$\psi(x) := \phi(x)h$ and $\psi_\#(\tau) = h^{-1} \phi_\#(\tau) h$, then $\psi \colon Q \to M$ and
is $\psi_\#$ equivariant.
We will call the (infinite dimensional) orbifold given by  the groupoid $[\left(\LL [M/G]\right)/G]$ the \emph{loop orbifold}.
\end{definition}

We need to consider the  equivalent over an orbifold of a Riemann surface with boundary.
This will consist of a map $\Phi : P \to M$ of a
$\Gamma$-principal bundle $P$ over an oriented Riemann surface $\Sigma$ ($\Gamma$ finite) and a homomorphism
$\Phi_\# \colon \Gamma \to G$ such that $\Phi$ is $\Phi_\#$-equivariant. Note that there is a natural
action of the group $G$ on $\Phi$. It is defined in the same way as for loops.

To define string connections in the case of orbifolds we must deal in one way or with 2-categories. Roughly speaking we \emph{define} $\SSS^1(\Xx)$ as a 2-category where the objects are orbifold loops $(\phi,\phi_\#)$, the arrows are orbifold surface maps as above. Then
the boundary $\partial P$ of $P$ will consist of $p$ incoming orbifold loops $\gamma_i \colon Q_i \to M$ $1\leq i \leq p$
with the induced orientation,
and $q$ outgoing ones $\gamma_j \colon \overline{Q}_j \to M$, $p+1 \leq j \leq p+q$
with the opposite orientation so that $\partial P = \bigsqcup_i Q_i \sqcup \bigsqcup_j Q_j$.
Here the  $Q_i$'s and the $Q_j$'s are $\Gamma$-principal bundles
over the circle. The 2-morphism of the 2-category  are given by the natural action of $G$ on the orbifold surface maps.
 We will define an orbifold string connection for $\Xx=[M/G]$ to be a 2-functor $\SSS^1(\Xx) \To \VS_1(\complex)$, namely a $G$-equivariant ordinary functor.

In \cite{LupercioUribeHolonomy} we prove the following refined version of the transgression (for a general orbifold $\Gg$).

\begin{theorem}
There is a natural holonomy homomorphism
$$\tau_2 :\breve{C}^2(\Gg,\U{1}(3)) \To \breve{C}^1(\Loop \Gg,\U{1}(2))$$
from the group of gerbes with connection over the orbifold $\Gg$ to the group of
line bundles with connection over the loop groupoid. Moreover this
  holonomy map commutes with the coboundary operator
and therefore induces a map in orbifold Beilinson-Deligne cohomology
$$H^3(\Gg;\integer(3))
\To H^2(\Gg;\integer(2)).$$
\end{theorem}

In fact we give a proof for the corresponding statement in $n$-gerbes. So given a gerbe $L=(g,A,B)$ we obtain a line orbibundle $E$ over the loop orbifold $\Loop \Xx$.

\begin{definition} \label{def.inertiagroupoid}
The inertia groupoid $\wedge \Gg$ is defined by:
\begin{itemize}
\item Objects $(\wedge \Gg)_0$: Elements $v \in \Gg_1$ such that $s(v) = t(v)$.
\item Morphisms $(\wedge \Gg)_1$: For  $v,w \in (\wedge \Gg)_0$ an
arrow $v \stackrel{\alpha}{\to} w$ is an element $\alpha \in
\Gg_1$ such that $v \cdot \alpha = \alpha \cdot w$
        $$
         \xymatrix{
         \circ \ar@(ul,dl)[]|{v} \ar@/^/[rr]|{\alpha}
         &&\circ \ar@(dr,ur)[]|{w^{-1}} \ar@/^/[ll]|{\alpha^{-1}}
         }$$
\end{itemize}
\end{definition}

One of the main
results of \cite{LupercioUribeConfiguration} is the following theorem

\begin{theorem} The fixed suborbifold of $\Loop \Gg$ under the
natural $S^1$-action (rotating the loops) is $$ \wedge \Gg = (\Loop
\Gg)^{S^1}$$
\end{theorem}

The following definition is due to Ruan \cite{Ruan1, Ruan3, Ruan5}. He used this definition to obtain  a twisted version of the Chen-Ruan cohomology \cite{ChenRuan} that has revived the interest in the theory of orbifolds in the last few years.

\begin{definition}
An ``inner local system'' is a flat line bundle
$\LL$ over the inertia groupoid $\wedge \Gg$ such that:
\begin{itemize}
\item $\LL$ is trivial once restricted to $\ident(\Gg_0) \subset \wedge \Gg_1$ (i.e. $\LL|_{\ident(\Gg_0)} = 1)$ and
\item $i^* \LL = \LL^{-1}$ where $i : \wedge \Gg \to \wedge \Gg$ is the inverse map (i.e. $(i(v,\alpha) = (\alpha^{-1} v \alpha, \alpha^{-1})$).
\end{itemize}
\end{definition}

\begin{theorem} The restriction of the holonomy of a gerbe with
connection over $\wedge \Gg$ (that is a line bundle with connection over
$\Loop \Gg$) is an inner local system on $\wedge \Gg$.
\end{theorem}

In the case in which we have a Lie group acting with finite stabilizers these line bundles are
the coefficients Freed-Hopkins-Teleman \cite{FreedHopkinsTeleman}
used to twist the cohomology of the twisted sectors in order to
get a Chern character isomorphism with the twisted $K$-theory of the
orbifold. We have used gerbes in \cite{LupercioUribeKtheory} to obtain twisted versions of $K$-theory that act a recipients of the charges of $D$-branes in string theory \cite{WittenKtheory} generalizing work of Adem and Ruan
\cite{AdemRuan}.

Returning to the subject of string connections we have the following result.

\begin{theorem}
Take a global gerbe $\xi$ with connection over $[M/G]$ and let $E$ be  the line
bundle with connection induced by it via transgression. Then $\xi$ permits to define
a string connection $\hol$ extending the line bundle $E$ of the loop groupoid $(\LL [M/G])/G$.
\end{theorem}

The analogous result for a general orbifold is more subtle and we refer the reader to \cite{LupercioUribeStringII} for details. There we use this theorem to generalize the results of Freed and Witten \cite{FreedWitten} on anomaly cancellation in string theory to the orbifold case.

To conclude let us mention that building on an idea of Hopkins and Singer \cite{HopkinsSinger} we have defined orbifold
Chern-Simons cohomology. The main difficulty here is to make sense of what an orbifold differential character should be
 \cite{LupercioUribeStringI}. We make a definition in such a way that we can prove the following result.

\begin{theorem}
The orbifold Beilinson-Deligne cohomology and the orbifold Cheeger-Simons
cohomology are canonically isomorphic.
\end{theorem}

\section{Acknowledgments}

This paper is based on lectures given by the first author at the IMATE-UNAM-Cuernavaca
and by the second author at the summer conference on
Geometric and topological methods for quantum field theory held at Villa de Leyva, Colombia.
The first author would like to thank the invitation and hospitality of Jos\'{e} Seade to the IMATE. The second
author would like to thank the organizers of the summer conference for the invitation
and the hospitality in Colombia, and especially to Sylvie Paycha for making a reality
 the publication of the present volume.
Both authors would also like to thank conversations with A. Adem, L. Borisov,
I. Kriz, I. Moerdijk, T. Nevins, M, Poddar, Y. Ruan and G. Segal.

\bibliographystyle{plain}
\bibliography{LupercioUribe_ambp}
\end{document}